%% file: RBF_GL.tex
\documentclass[a4paper,12pt]{article}
\input{Util.tex}

\usepackage{authblk}

\usepackage{etoolbox}

\patchcmd{\thebibliography}{%
  \section*{\refname}%
}{%
  \subsection*{\refname}%
}{}{}

\apptocmd{\thebibliography}{%
  \setlength{\itemsep}{0pt}%
  \setlength{\parskip}{0pt}%
}{}{}

\graphicspath{ {Figures/} }

\numberwithin{equation}{section}

\DeclareMathOperator{\tr}{tr}
\title{Anisotropic Radial Basis Function Methods for Continental Size Ice Sheet Simulations}

\author{Gong Cheng\footnote{Uppsala University, Division of Scientific Computing, cheng.gong@it.uu.se} $\,$ and Victor Shcherbakov\footnote{Uppsala University, Division of Scientific Computing, victor.shcherbakov@it.uu.se} }

\begin{document}
\maketitle

\begin{abstract}
In this paper we develop and implement anisotropic radial basis function methods for simulating the dynamics of ice sheets and glaciers. We test the methods on two problems: the well-known benchmark ISMIP-HOM B that corresponds to a glacier size ice and a synthetic ice sheet whose geometry is inspired by the EISMINT benchmark that corresponds to a continental size ice sheet. We illustrate the advantages of the radial basis function methods over a standard finite element method. We also show how the use of anisotropic radial basis functions allows for accurate simulation of the velocities on a large ice sheet, which was not possible with standard isotropic radial basis function methods due to a large aspect ratio between the ice length and the ice thickness. Additionally, we implement a partition of unity method in order to improve the computational efficiency of the radial basis function methods.
\end{abstract}

\section{Introduction}

There is a growing interest in simulating the evolution of ice sheets for predicting their contribution to the future sea level rise \cite{pachauri2014climate} and also understanding the processes of forming past landscapes. 
Mathematical models are introduced as tools to understand the dynamics of ice sheets in the past and in the future \cite{deconto2016contribution}.
 Improvement in accuracy and efficiency of the modeling and numerical methods is always needed, especially for large scale and long time simulations \cite{vaughan2007hard, pattyn2008benchmark, pattyn_results_2012}.

The ice flow is generally described as an incompressible, non-Newtonian fluid with highly nonlinear viscosity.  
One accurate model for simulating the ice dynamics is the so-called full Stokes equations \cite{greve2009dynamics, baral2001asymptotic, hutter1983theoretical}. 
The deformation of the ice body under its own weight is governed by Glen's flow law \cite{glen1958flow}. It relates the stress field to the strain rates as a viscous fluid and the viscosity depends nonlinearly on the velocities, which introduces an additional degree of difficulty in the numerical solving procedure for the full-Stokes equations. Moreover, the discretisation of the full-Stokes system gives rise to a saddle-point problem, which in a setting for finite element methods requires special numerical treatments for the methods to satisfy the inf-sup condition \cite{larson2013finite}, such as adding stabilization on the pressure variables \cite{blasco2007anisotropic} or using high-order finite element methods \cite{isaac2015solution}.
Therefore, the numerical solution of the full-Stokes equations is demanding in terms of computational time.

Several simplifications are derived for the full-Stokes equations to reduce the computational complexity. 
The first order Stokes model (also known as the Blatter--Pattyn model) is based on the assumption that the hydrostatic pressure is balanced by the vertical normal stress \cite{Blatter95,Pattyn2003}, such that the horizontal gradient of the vertical velocity is neglected. The system is simplified to an elliptic problem that only contains the horizontal velocities as unknowns and the vertical velocity is recovered by solving the incompressibility equation.  
Other approximations are the Shallow Ice Approximation (SIA) \cite{BlaHaftet} and Shallow Shelf Approximation (SSA) \cite{SSA,WGH99}, which however give a lower order of approximation than the first order Stokes model, i.e., have a larger model error.

These models are intercompared within several benchmark experiments during the past decade. For instance, some of the well-known benchmark experiments are the Ice Sheet Model Intercomparison Project for Higher-Order ice sheet Models (ISMIP-HOM) \cite{pattyn2008benchmark} and the framework of European Ice Sheet Modeling Initiative (EISMINT) \cite{huybrechts_eismint_1996}. In ISMIP-HOM, different Stokes approximations are compared on glacier size problems (about 10~km long), whereas in EISMINT, the computational domains are continental size (more than 1000~km long). 

Traditional numerical methods such as finite element methods (FEM) are commonly used for solving ice sheet models since FEM can easily handle complex geometry with different types of boundary conditions. However, it has some drawbacks when solving problems on a domain with a moving boundary, which leads to remeshing of the entire domain for every time step. Also, solving the nonlinear system requires a full matrix reassembly during each nonlinear iteration. Therefore, a mesh-free method that can be stated in strong form is considered to be a preferred choice for such problems. Radial basis function (RBF) methods are of that kind~\cite{Fasshauer}. The idea is to define a finite-dimensional basis, which consists of functions, whose values depend on the distance from their centers, and use them for approximating the solution.

RBF methods have been introduced for ice simulations in \cite{Ahlkrona}, where the authors have studied the advantages of an RBF approach for the Haut glacier d'Arolla, which is also a test case from the ISMIP-HOM benchmark. In this paper, we continue the work and extend the approach to solve for dynamics of continental size ice sheets. We introduce  anisotropic RBF approximations to solve problems with continental ice sheet geometries, which typically have large aspect ratios. The high aspect ratio is an obstacle for standard isotropic RBFs due to the strong dependence of the shape parameter, which determines the width of the functions. The use of anisotropic RBFs significantly relaxes this dependence and simplifies the method implementation. Additionally, we suggest a strategy of selecting the shape parameter value based on the conditioning of the interpolation matrix. In order to enable more efficient simulations we develop a partition of unity approach that is also based on anisotropic patches. Also, we study error estimates for the anisotropic method and perform convergence tests.


The remainder of the paper is structured as follows. 
In Section~\ref{sec:ice} we explain the full-Stokes equations and the first order Stokes model that governs the dynamics of the ice sheets. Then, in Section~\ref{sec:rbf} we introduce the anisotropic RBF methods, and in Section~\ref{sec: results}, we present the numerical results obtained by the RBF methods for the two test cases. We also provide a comparison of the RBF methods with the standard FEM in terms of time-to-accuracy, i.e., we compare the execution times to achieve certain error levels.
Finally, we draw some conclusions in Section~\ref{sec:conclusion}.

\section{Ice Sheet Dynamics}\label{sec:ice}
Ice can be viewed as a very viscous fluid flow on a large scale. Thus, the models for simulating ice sheet dynamics are inspired by fluid dynamics laws.
In general, ice is considered as an incompressible flow with a low Reynolds number and with the stress tensor related to the strain rate by a power law viscous rheology \cite{greve2009dynamics}. Due to the slow motion of ice masses, the acceleration term can be neglected and the Navier--Stokes equations can be turned into the full-Stokes equations that describes a steady flow. However, the computational demand for solving the full-Stokes system can be excessive, especially for such large domains as ice sheets. Therefore, under assumptions of large aspect ratios between the ice length and ice thickness, the horizontal variations of the vertical velocity in the full-Stokes equations can be neglected. Thereby, we arrive at the simplified first order Stokes equations, which is an approximation to the full-Stokes equations in terms of the thickness/length ratio. The solution obtained by the first order Stokes equations is accurate at the ice margins, but the areas with discontinuities in the surface gradient might not be resolved well enough, if compared with the solution of the full-Stokes system.

\begin{figure}[!ht]
\centering
\includegraphics[width=0.9\textwidth]{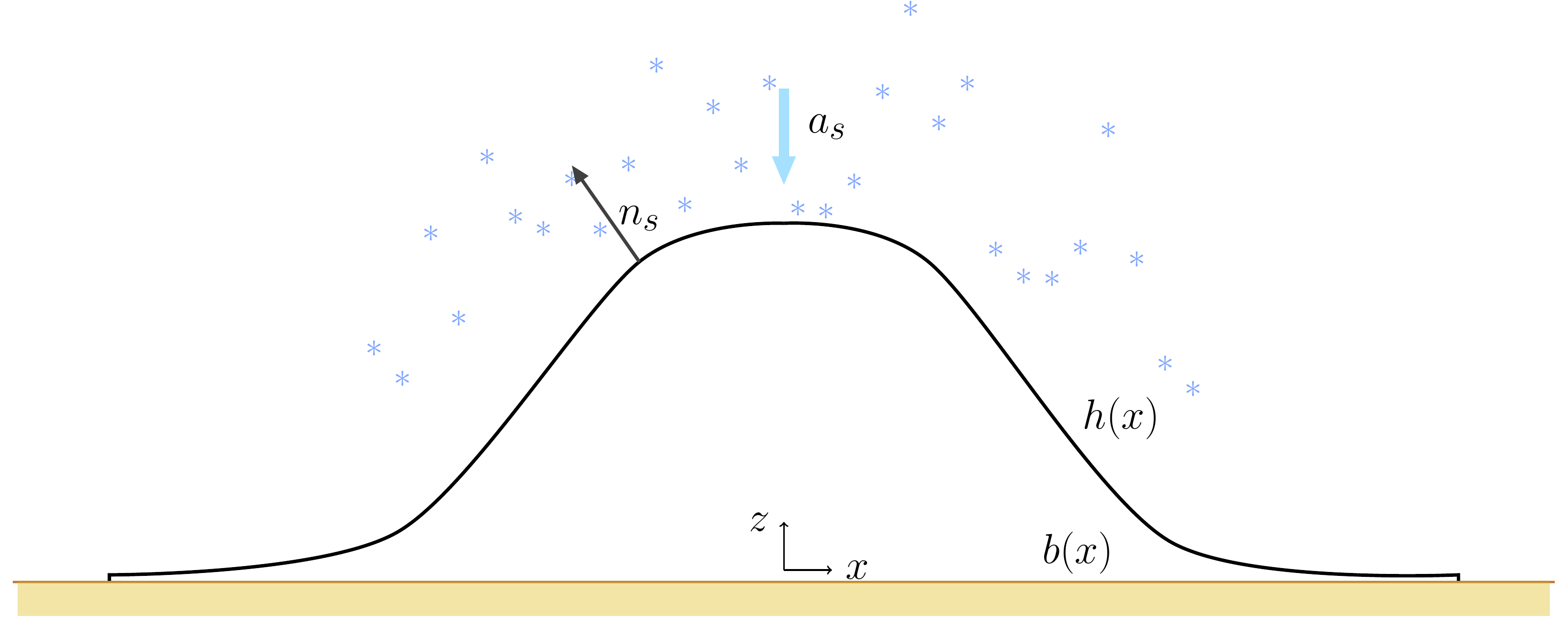}\\
\caption{A schematic continental size ice cap. For aesthetic reasons the coordinates are exaggerated in the vertical direction.}
\label{Icesheet}
\end{figure}


\subsection{The Full-Stokes Equations}\label{sec:Stokes}
The nonlinear full-Stokes equations are defined by the conservation of momentum and mass
\begin{equation}
\begin{array}{r@{\quad }c@{\quad }l}
\nabla \cdot \big(\eta (\nabla \mathbf{v}+\nabla \mathbf{v}^{T})\big)-\nabla
p+\rho \mathbf{g}=0,
\\[5pt]
\nabla \cdot \mathbf{v}=0, 
\end{array}
\label{eq:FS}
\end{equation}
where $\mathbf{v}$ is the vector of velocities $\mathbf{v}=\left(
\begin{matrix}
v_{x}&v_{y}&v_{z}\end{matrix}\right)  ^{T}$, $\rho $ is the density of the ice, $\eta$ is the viscosity, $p$ is the pressure, and $\mathbf{g}$ is the gravitational acceleration.
A constitutive equation (so-called Glen's flow law) relates the deviatoric stress tensor $\mathbf{{T}}^{D}$ and the strain rate $\mathbf{D}=\frac{1}{2}\left(  \nabla \mathbf{v}+(\nabla\mathbf{v})^{T}\right)$ such that
\begin{equation}
  {\mathbf{{D}}}=\mathcal{A}(T^{\prime })f(\sigma )\mathbf{{T}}^{D},
  \label{eq:glen}
\end{equation}
where $\mathcal{A}(T^{\prime })$ is the rate factor that describes how the viscosity depends on the pressure melting point corrected temperature $T^{\prime }$. For isothermal flow, the rate factor $\mathcal{A}$ is constant.
The value of the physical parameters are given in Table \ref{tab: physical param}.

The deviatoric stress tensor $\mathbf{{T}}^{D}$ is given by
\begin{equation}
{\mathbf{{T}}}^{D}=\left(
\begin{array}{c@{\quad }c@{\quad} c}
t_{xx}^{D}&t_{xy}^{D}&t_{xz}^{D}
\\[2pt]
t_{yx}^{D}&t_{yy}^{D}&t_{yz}^{D}
\\[2pt]
t_{zx}^{D}&t_{zy}^{D}&t_{zz}^{D}
\end{array}
\right)  ,
\label{eq:stress}
\end{equation}
where $t_{xx}^{D}$, $t_{yy}^{D}$, $t_{zz}^{D}$ and $t_{xy}^{D}$ denote longitudinal stresses and $t_{xz}^{D}$, $t_{yz}^{D}$ vertical shear stresses. The function $f(\sigma)$ is called the creep response function for ice with
\begin{equation}
f(\sigma )=\sigma^{n-1},
\label{eq:creep}
\end{equation}
where $n$ is the flow-law exponent factor and $\sigma $ is the effective stress, which is defined as the second invariant of the deviatoric stress tensor
\begin{equation}
\sigma=\left[\big(t_{x{y}}^{D}\big)^{2}+\big(t_{yz}^{D}\big)^{2}+\big(t_{xz}^{D}\big)^{2}+
\frac{1}{2}\Big(
\big(t_{xx}^{D}\big)^{2}+\big(t_{yy}^{D}\big)^{2}+\big(t_{zz}^{D}\big)^{2}\Big)\right]^{\frac{1}{2}}  .
\label{eq:effectivestress}
\end{equation}
The effective strain rate is defined in analogy with the effective stress such that
\begin{equation}
\mathbf{D}_e= \sqrt{\frac{1}{2}\tr\big(\mathbf{D}\cdot\mathbf{D}\big)},
\label{eq:effectivestrain}
\end{equation}
with the relation
\begin{equation}
  \mathbf{D}_e = \mathcal{A}\sigma^n,
\end{equation}
and the viscosity $\eta $ is defined by
\begin{equation}
\eta =\big( 2\mathcal{A}f(\sigma )\big)  ^{-1}=\frac{1}{2}\mathcal{A}^{-\frac{1}{n}}\mathbf{D}_e^{\frac{1-n}{n}}.
\label{eq:eta}
\end{equation}
%
At the top surface with normal {{$\mathbf{n}_{s}$}}, the ice is stress-free, i.e.,
\begin{equation}
\big(-p\mathbf{I}+2\eta\mathbf{D}\big)\cdot {\mathbf{n}_{s}}=0,
\label{eq:stressfree}
\end{equation}
where $\mathbf{I}$ is the identity matrix.

At the base of the ice, the normal vector {$\mathbf{n}_{b}$} and the tangential vectors $\mathbf{t}_{1}$ and $\mathbf{t}_{2}$ span the base surface such that $\mathbf{n}_{b}\cdot \mathbf{t}_{i} = 0$, $i=1,2$ and $\mathbf{t}_{1} \cdot \mathbf{t}_{2}=0$.  If the ice base is frozen, the velocity $\mathbf{v}$ satisfies a no slip condition at the base
\begin{equation}
\mathbf{v}=0.
\label{eq:noslip}
\end{equation}
%


\begin{table}[h!]
  \begin{center}
    \begin{tabular}{ c c c l }
      \hline\hline
      Parameter & Value & Unit &  Description \\
      \hline
      $\rho$ & 900 & kg m$^{-3}$ & Ice density\\
      $g$ & $9.81$ & m s$^{-2}$  & Acceleration of gravity \\
      $n$ & $3$ & -- & Flow-law exponent \\
      $\mathcal{A}$ & $10^{-16}$ & Pa$^{-3}$a$^{-1}$ & Rate factor in flow law \\
      \hline\hline
    \end{tabular}
    \caption{The physical parameters of ice.} \label{tab: physical param}
  \end{center}
\end{table}

The time evolution of the domain geometry is introduced by a kinematic boundary condition at the top surface, which moves the surface $h(x,y)$ according to the surface velocities and the accumulation-ablation function $a_s$ such that
\begin{equation}
  \pdd{h}{t}+v_x\pdd{h}{x}+v_y\pdd{h}{y}=a_s+v_z,
\end{equation}
where $a_s$ is in the unit of meters per annum (m/a) ice equivalent. Positive values of $a_s$ imply snowing and negative values imply snow melting.

\subsection{The First Order Stokes Model}

As the full-Stokes system consists of complex equations with a highly nonlinear viscosity, they require a massive computational effort to be solved numerically. Moreover, the system gives rise to a saddle point problem that adds an extra degree of difficulty. Therefore, the full-Stokes equations are often simplified to a reduced form under the assumption that the variational stress is neglected \cite{pattyn2003new} due to the large aspect ratio. Thus, the conservation of momentum in \eqref{eq:FS} becomes
\begin{equation}
  \begin{split}
    \pdd{}{x}\left(4\eta\pdd{v_x}{x}+2\eta\pdd{v_y}{y}\right)+\pdd{}{y}\left(\eta\pdd{v_x}{y}+\eta\pdd{v_y}{x}\right)+ \pdd{}{z}\left(\eta\pdd{v_x}{z}\right)&=\rho g\pdd{h}{x},\\
    \pdd{}{x}\left(\eta\pdd{v_x}{y}+\eta\pdd{v_y}{x}\right)+\pdd{}{y}\left(4\eta\pdd{v_y}{y}+2\eta\pdd{v_x}{x}\right)+ \pdd{}{z}\left(\eta\pdd{v_y}{z}\right)&=\rho g\pdd{h}{y},\\
  \end{split}\label{eq:first order stokes}
\end{equation}
where $h(x,y)$ denotes the surface elevation and the viscosity is written as 
\begin{equation}
  \begin{split}
    \eta = &\frac{1}{2}\mathcal{A}^{-\frac{1}{n}}\left[\left(\pdd{v_x}{x}\right)^2+\left(\pdd{v_y}{y}\right)^2+\pdd{v_x}{x}\pdd{v_y}{y}\right.\\
    &+\left.\frac{1}{4}\left(\pdd{v_x}{y}+\pdd{v_y}{x}\right)^2+\frac{1}{4}\left(\pdd{v_x}{z}\right)^2+\frac{1}{4}\left(\pdd{v_y}{z}\right)^2 \right]^{\frac{1-n}{2n}}.
  \end{split}
\end{equation}
The vertical velocity is obtained through vertical integration over the conservation of mass in \eqref{eq:FS} from the bottom $b(x,y)$ to a height $z(x,y)$
\begin{equation}
  v_z(z)-v_z(b)=-\int_{b}^{z} \left(\pdd{v_x}{x}+\pdd{v_y}{y}\right)\text{d} \xi.
  \label{eq:integralZ}
\end{equation}
The stress-free boundary on the top surface is expressed in terms of the velocity gradients such that
\begin{equation}
  \begin{split}
    \left(4\eta\pdd{v_x}{x}+2\eta\pdd{v_y}{y}\right)\pdd{h}{x}+\left(\eta\pdd{v_x}{y}+\eta\pdd{v_y}{x}\right)\pdd{h}{y}-\eta\pdd{v_x}{z}&=0,\\
    \left(4\eta\pdd{v_y}{y}+2\eta\pdd{v_x}{x}\right)\pdd{h}{y}+\left(\eta\pdd{v_x}{y}+\eta\pdd{v_y}{x}\right)\pdd{h}{x}-\eta\pdd{v_y}{z}&=0.
  \end{split}\label{eq: stress free 3d}
\end{equation}


The system of equations that has to be solved is referred to as the first order Stokes or the Blatter--Pattyn model. It is second order accurate with respect to the thickness/length ratio.  The computational demand for solving the Blatter--Pattyn model is significantly lower than for solving the full-Stokes system. Since in this paper we are not dealing with simulations of grounding lines and calving fronts, we use the Blatter--Pattyn model which  gives a satisfactory solution.

\subsection{The Flow-Line Model}

We consider the two-dimensional flow-line model, which is a vertical cutting plane along the surface gradient. This is a simplification from the three-dimensional ice sheet to a two-dimensional problem, which is commonly used in glaciology. The Blatter--Pattyn equations in \eqref{eq:first order stokes} combined with the conservation of mass in \eqref{eq:FS} are reduced to 
\begin{equation} \label{eq: B-P model}
  \begin{split}
  4\pdd{}{x}\left(\eta \pdd{v_x}{x}\right)+  \pdd{}{z}\left(\eta \pdd{v_x}{z}\right)&=\rho g\pdd{h}{x},\\
  \pdd{v_z}{z} &= -\pdd{v_x}{x},
  \end{split}
\end{equation}
with the viscosity (using $n=3$ in Table \ref{tab: physical param})
\begin{equation} \label{eq: viscosity}
  \eta=\frac{1}{2}\mathcal{A}^{-1/3}\left[\left(\pdd{v_x}{x}\right)^2+\frac{1}{4}\left(\pdd{v_x}{z}\right)^2\right]^{-1/3}.
\end{equation}
The horizontal velocity $v_x$ is  determined by the conservation of momentum in \eqref{eq: B-P model} and the vertical velocity $v_z$ can be computed by vertically integrating~$v_x$ from the bottom of the ice to the  position $z$ as in \eqref{eq:integralZ} 
\begin{equation}
  v_z(z)-v_z(b)=-\int_{b}^{z} \pdd{v_x}{x} \text{d} z.
\end{equation}

The geometry of the ice bottom is constrained by the bedrock and the bedrock is assumed to have no deformation during the whole simulation. There is no sliding on the bedrock, therefore the horizontal and vertical velocity at $z=b$ are both equal to zero.
%
On the top surface, the stress-free boundary condition is simplified as 
\begin{equation} \label{eq: top bc}
  \eta\left(4\pdd{v_x}{x}\pdd{h}{x}-\pdd{v_x}{z} \right) = 0.
\end{equation}
The lateral boundary conditions are given in the following ways for different test cases: 
\begin{enumerate}

  \item For the periodic boundary conditions, e.g., in the ISMIP-HOM benchmark experiment, the solution from the right boundary is mapped to the left boundary.

  \item At the far end of the thin ice area, e.g., the left and right boundaries in the ice cap experiment, a symmetric boundary condition is imposed as 
  \begin{equation} 
    \pdd{h}{x}=0, \quad v_x=0.
  \end{equation}
\end{enumerate}

\section{Radial Basis Function Methods}\label{sec:rbf}

In order to solve the system~\eqref{eq: B-P model} numerically, we develop and implement RBF methods. To construct an RBF approximation we scatter a set of nodes over the computational domain, where each node is associated with a basis function. The collection of basis functions forms a finite basis in the functional space. Some typical choices of basis function can be found in Table~\ref{TabRBF}. An important feature of an RBF is that it depends only on the distance between the nodes and its center. This is a valuable property since it makes the method  easily applicable to high dimensional problems. Additionally, RBF methods are mesh-free and therefore suitable for problems which are defined in domains with complex geometries, such as ice sheets and glaciers.  
\begin{table}[h]
\begin{center}
\caption{Commonly used smooth radial basis functions.}
\label{TabRBF}
\begin{tabular}{ l  c  c  c  l  }
\hline\hline 
RBF & & &  & $\phi(\varepsilon,r)$ \vspace{5pt}  \\ \hline 
 Multiquadric (MQ) &  & &  & $(1+(\varepsilon r)^2)^{1/2}$ \\ 
 Inverse Multiquadric (IMQ) & & &  & $(1+(\varepsilon r)^2)^{-1/2}$ \\
 Inverse Quadratic (IQ) & & &  & $(1+(\varepsilon r)^2)^{-1}$ \\
 Gaussian (GA) &  & &  &  $e^{-(\varepsilon r)^2}$ \\
\hline\hline
\end{tabular}
\end{center}
\end{table}

Given $N$ distinct scattered nodes $\underline{x}=\{x_{1},\ldots,x_{N}\}$, $x_{i}\in \Omega\subset \mathbb{R}^{d}$, we can construct an RBF approximation $\tilde v(x)$ of a function $v(x)$ with values $\underline{v}=[v(x_{1}),\ldots,v(x_{N})]$ defined at the nodes such that
\begin{equation}\label{RBFinterp}
\tilde v(x) = \sum_{j=1}^{N} \lambda_{j}\phi(\varepsilon,||x-x_{j}||), \quad x\in\Omega,
\end{equation}
where $\lambda_{j}$ are the unknown coefficients, $\|\cdot\|$ is the Euclidean norm and $\phi(\varepsilon,r)$ is a real-valued radial basis function, whose width is determined by the shape parameter $\varepsilon$. In order to determine the coefficients $\lambda_{j}$ we collocate the approximation~(\ref{RBFinterp}) at the node set and obtain a system of linear equations for the coefficients $\underline{\lambda}$
\begin{equation}\label{systemA}
A\underline{\lambda} = \underline{v},
\end{equation}
where the matrix $A$ has  the following elements ${A_{ij} = \phi(\varepsilon,\|x_i-x_j\|)}$.

\subsection{Anisotropic Radial Basis Functions}

\begin{figure}[!ht]
\centering
\includegraphics[width=0.9\textwidth]{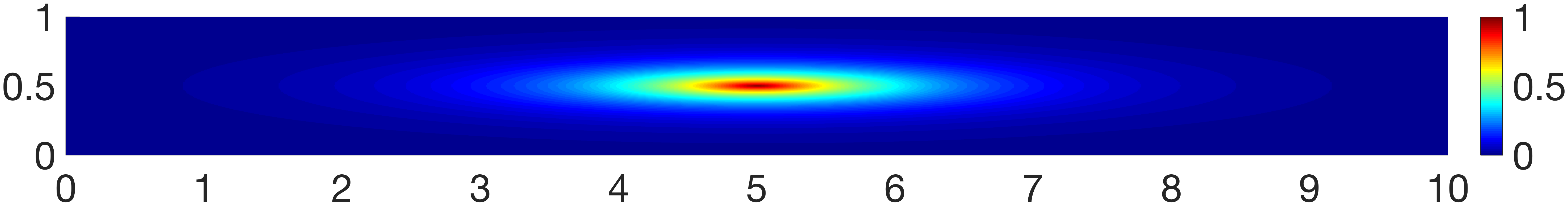}\\
\vspace{0.3cm}
\includegraphics[width=0.9\textwidth]{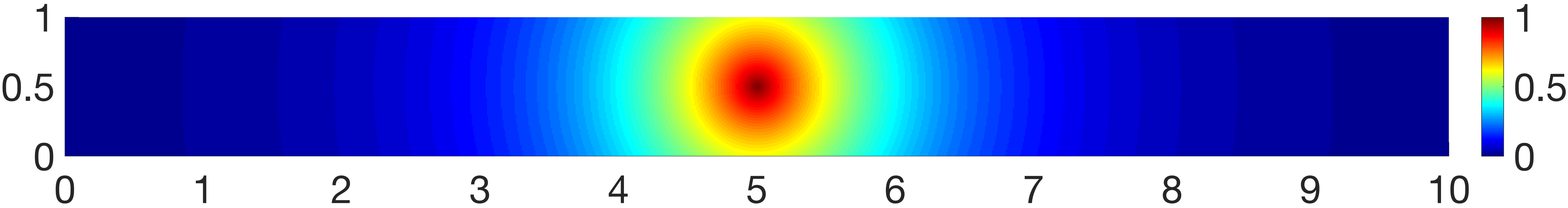}
\caption{\emph{Top:} An anisotropic Gaussian RBF. \emph{Bottom:} An isotropic Gaussian RBF. The aspect ratio $a=10$.}
\label{FlatRBF}
\end{figure}

The thickness of a continental ice sheet is relatively small in comparison with its length. The  ratio between the length and thickness may in some cases reach $500:1$. Therefore, approximation with standard RBFs fails to provide a reasonable resolution in the vertical direction. Several authors suggested to use anisotropic (elliptic) RBFs instead~\cite{Beatson2010,Bayona,Zhao}. These are functions, which are scaled in an appropriate way to give good resolution and match the domain geometry features (see Figure~\ref{FlatRBF}). This can be implemented by redefining the distance between two nodes such that
\begin{equation}\label{dist}
\|x-y\|_{a} = \sqrt{a_1^2(x_1-y_1)^2 + \cdots  + a_d^2(x_d-y_d)^2}, \quad x, y \in \Omega \subset \mathbb{R}^d,
\end{equation}
where $a_i$ are the aspect ratios of the domain discretisation defined as the ratio between typical node distances in the respective directions (see Section~\ref{sec: results} for an accurate definition).
Without loss of generality, we assume  $a_1=1$  and $a_i\geq 1$, $i=2,\ldots,d$. Particularly, for the two-dimensional problems that we consider in this paper the vertical dimension is several orders of magnitude smaller than the horizontal dimension. Therefore, the norm can be defined as 
\begin{equation}
\|x-y\|_{a} = \sqrt{(x_1-y_1)^2 + a^2(x_2-y_2)^2}, \quad x,y \in \Omega \subset \mathbb{R}^2. \label{eq: a_norm}
\end{equation}
Thus, the basis function $\phi$ is no longer a function of the Euclidean distance but a function of the scaled distance defined by the $\|\cdot\|_a$-norm. 

%
%

%
%

\subsection{Kansa's Method}
The RBF interpolant~(\ref{RBFinterp}) can also be collocated for a partial differential equation (PDE) or a system of partial differential equations. This approach is known as Kansa's method~\cite{Kansa1,Kansa2}. The nonlinear First Order Stokes Equations~\eqref{eq: B-P model} can be written in general form as a nonlinear boundary value problem
\begin{equation}\label{nonBVP}
\mathcal{P}\left[x,v(x),\mathcal{D}v(x)\right] = 0 \quad \Longleftrightarrow \quad
\begin{cases}
\mathcal{P}_1 = 0, & x \in \Omega,\\
\mathcal{P}_2 = 0, & x \in \partial \Omega,
\end{cases}
\end{equation}
where $\mathcal{P}_1$ is the interior nonlinear operator,  $\mathcal{P}_2$ is the boundary nonlinear operator, and $\mathcal{D}$ is a shorthand notation for differential operators, such as $\partial_x$, $\partial_z$, $\nabla$. 

A root of the nonlinear system~\eqref{nonBVP} can be sought by a nonlinear solver, such as Newton's method~\cite{FasshauerNonLin} or a fixed point iteration method~\cite{Ahlkrona}, that iteratively solves a linearised problem. Thus, we arrive at a system of linear equations in the following form
\begin{equation}\label{BVP}
\begin{cases}
\mathcal{L}v(x) = f(x), &x \in \Omega, \\
\mathcal{F}v(x) = g(x), &x \in \partial\Omega,
\end{cases}
\end{equation}
where $\mathcal{L}$ is the linearised interior differential operator, $\mathcal{F}$ is the linearised boundary differential operator, and $f$, $g$ are the right hand side forcing functions. We seek a solution to system~(\ref{BVP}) in the form of the RBF~interpolant~(\ref{RBFinterp}). 
Collocating at the node points we obtain the following system of linear equations
\begin{equation}\label{discreteBVP}
C \underline{\lambda} := 
\begin{bmatrix}
L\\
F
\end{bmatrix}
\underline \lambda := 
\begin{bmatrix}
L_{II} & L_{IB}\\
F_{BI} & F_{BB} 
\end{bmatrix}
\begin{bmatrix}
\underline \lambda_{I} \\
\underline \lambda_{B}
\end{bmatrix} = 
\begin{bmatrix}
\underline f_{I} \\
\underline g_{B}
\end{bmatrix},
\end{equation}
where $L$, $F$, $\underline f$, and $\underline g$ are discrete representations of the continuous quantities, and the subscripts $I$ and $B$ denote that the quantities are evaluated on the interior and the boundary nodes, respectively (without loss of generality, we assume that the first $N_{I}$ nodes belong to the interior of $\Omega$ and the last $N_{B}=N-N_{I}$ nodes belong to the boundary $\partial \Omega$). The matrices $L$ and $F$ are constructed from the elements $L_{ij} = \mathcal{L}\phi(\varepsilon,\|x_i-x_j\|_a)$ and $F_{ij} = \mathcal{F}\phi(\varepsilon,\|x_i-x_j\|_a)$. 

It has been shown~\cite{Driscoll,Larsson7,Schaback2} that for smooth RBFs the magnitude of the coefficients becomes unbounded as $\varepsilon \to 0$, while the values $\tilde v(\underline{x})$ remain well-behaved. Therefore, we prefer to transform the problem into a search for the nodal values $\tilde v(\underline x)$. For the basis functions presented in Table~\ref{TabRBF} the interpolation matrix $A$ is non-singular for distinct node points~\cite{Micchelli}. Hence, based on \eqref{systemA} we can write
\begin{equation}\label{trans}
\underline\lambda = A^{-1}\underline{v}.
\end{equation}
Thus, relation~(\ref{trans}) allows us to transform the problem from solving for  $\underline \lambda$ to directly solving for $\underline v$. 
\begin{equation}\label{discreteBVPu}
\begin{bmatrix}
L_{II} & L_{IB}\\
F_{BI} & F_{BB} 
\end{bmatrix}
\begin{bmatrix}
\underline \lambda_{I} \\
\underline \lambda_{B}
\end{bmatrix} = 
\begin{bmatrix}
L_{II} & L_{IB}\\
F_{BI} & F_{BB} 
\end{bmatrix} A^{-1}
\begin{bmatrix}
\underline v_{I} \\
\underline v_{B}
\end{bmatrix}=
\begin{bmatrix}
\underline f_{I} \\
\underline g_{B}
\end{bmatrix}. 
\end{equation}
%
\subsection{Radial Basis Function Partition of Unity Method}
The approach presented in the previous section is referred to as a global RBF approximation, since it is constructed over all discretisation nodes. Such an approach gives a highly accurate approximation, but results in a dense system of linear equations, which is computationally expensive to solve. To overcome this issue we employ a partition of unity method  that allows for a significant sparsification of the linear system. Thereby, the high computational cost associated with the global method is reduced, while a similarly high accuracy is maintained~\cite{Ahlkrona,Shcherbakov}. Moreover, a partition based formulation is well suited for parallel implementations.

The partition of unity method was first introduced for finite element methods by Babu\v{s}ka and Melenk in~\cite{Babuska}. Later it was applied to RBF-based formulations by multiple authors~\cite{Fasshauer,Shcherbakov,Cavoretto,Safdari,Wendland2002}. The main idea of the method is to subdivide the computational domain into subdomains and construct an RBF interpolant locally in each subdomain and then combine them together by the partition of unity functions, which serve as weights. Below comes a formal description of the method.

To define a partition of unity method for problem~(\ref{BVP}) we construct a set of overlapping patches $\{\Omega_k\}_{k=1}^{M}$ that form an open \mbox{cover of the domain $\Omega$,} such that
\begin{equation}
\Omega \subset \bigcup_{k=1}^{M} \Omega_{k}.
\end{equation}
The patches $\Omega_k$ are selected as circular discs in the anisotropic norm $||\cdot||_a$. Additionally, we construct a partition of unity $\{w_k\}_{k=1}^{M}$ that is subordinated to the open cover $\{\Omega_k\}_{k=1}^{M}$. The function $w_k$ is compactly supported on $\Omega_k$, and
\begin{equation}\label{partun}
\sum_{k=1}^{M} w_k(x) = 1, \quad x \in \Omega.
\end{equation}
Thus, the RBF approximation can be written in the form of a weighted sum of all local approximations
\begin{equation}\label{LocToGlob}
\tilde v(x) = \sum_{k=1}^{M}w_k(x)\tilde v_k(x), \quad x \in \Omega,
\end{equation}
where $\tilde v_k(x)$ is a local approximation defined as
\begin{equation}\label{LocRBFinterp}
\tilde v_k(x) = \sum_{i = 1}^{n_{k}}\lambda_{i}^{k}\phi(\varepsilon,\|x-x_{i}^{k}\|_a), \quad x \in \Omega_k,
\end{equation}  
and $n_k$ is the local number of computational nodes in the patch $\Omega_k$.
The partition of unity weight functions can be constructed using Shepard's method~\cite{Shepard}
\begin{equation}\label{weightfunction}
w_k(x) = \frac{\varphi_{k}(x)}{\sum_{i=1}^{M}\varphi_{i}(x)}, \quad k = 1,\ldots,M.
\end{equation} 
In order to provide necessary smoothness of the solution we choose $\varphi_{k}(x)$ as a $C^2(\mathbb{R}^2)$ compactly supported Wendland function~\cite{Wendland}
\begin{equation}\label{Wendfunc}
\varphi(r) = 
\begin{cases}
  (1-r)^{4}(4r+1), & \text{if $0 \leq r \leq 1$},\\
  0, & \text{if $r>1$},
\end{cases}
\end{equation} 
which is scaled and shifted accordingly to fit the partition $\Omega_k$ with the centre $c_k$ and radius $\rho_k$
\begin{equation}
 \varphi_{k}(x) = \varphi_k \left( \frac{||x- c_k||_a}{\rho_k}\right), \quad \forall x \in \Omega.
\end{equation}


\begin{figure}[!ht]
  \center 
    \includegraphics[width=1\textwidth]{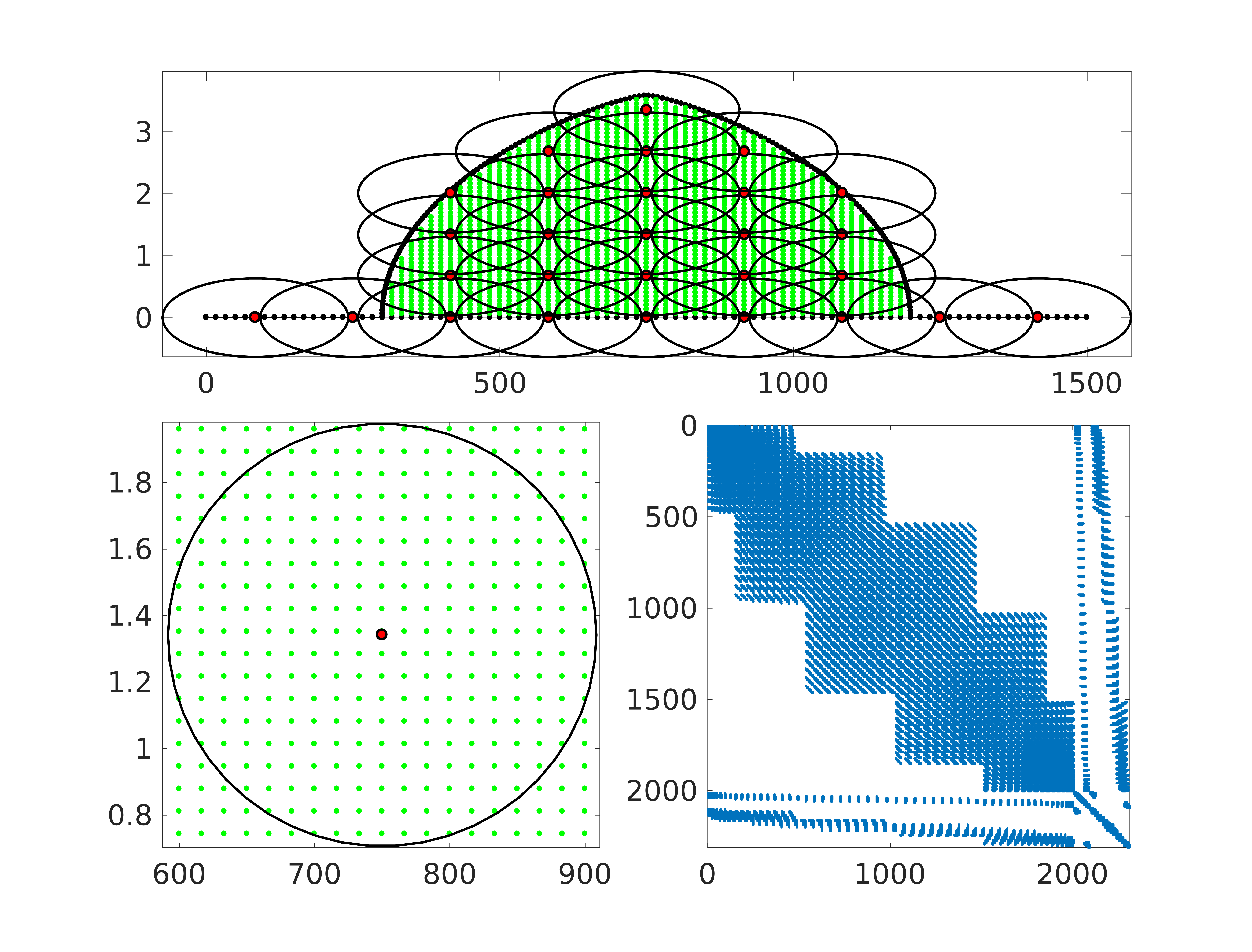}  
    \caption{\emph{Top:} An example of an anisotropic partitioning over the ice cap domain. \emph{Bottom left:} Zoom-in on one partition. \emph{Bottom right:} The sparse structure of the interpolation matrix for the ice cap domain with the discretization and partitioning as on the top panel.} \label{fig:EISMINT_PUM} 
\end{figure}

An example of the domain partitioning is shown in Figure~\ref{fig:EISMINT_PUM} with an ice cap geometry spanning $[0,1500]\times[0,3.5]$ kilometres. We use anisotropic partitions in this application as it was done in~\cite{Safdari}. The black ellipses are the boundaries of the partitions centered at the red dots and the green dots are the computational nodes in the background node set. 
A zoomed-in plot of one partition is shown on the bottom left panel in Figure~\ref{fig:EISMINT_PUM}. The ratio between the major and the minor axis of the partition is equivalent to the aspect ratio $a$. Consequently, the partitions can be considered as circularly shaped under the $||\cdot||_a$-norm given in \eqref{eq: a_norm}. 
On the bottom right panel, we show the sparsity structure of the interpolation matrix for the ice cap on $2311$ Cartesian nodes with the partitioning as in the top panel. 
There are $18.9\%$ nonzero element in the matrix with an average bandwidth around $800$ elements. The bandwidth will not grow with the number of computational nodes as long as the partition size is fixed to contain a certain number of nodes, e.g., for the finest grid that was used in our experiments for the ice cap test case the coefficient matrix had only $2.3\%$ of nonzero elements.




Note that the amount of overlap between patches is limited in order to maintain a sparse pattern of the linear system. Greater overlap leads to a slightly higher accuracy but denser linear system and, hence, higher computation costs. Less overlap leads to a lower accuracy but at the same time to a lower computational cost. We find that $20\% - 30\%$ overlap is appropriate. More on this matter can be found in~\cite{Shcherbakov}.


%

%
%
~
%
%

\subsection{Error Estimates}

We are interested in estimating the difference between the exact unknown solution $v$ and its RBF approximation $\tilde v$. We define the PDE error as
\begin{equation}\label{PDEerror}
E = v-\tilde v.
\end{equation}
We assume that we solve the already linearized problem \eqref{BVP} and the  error introduced by the linearization is negligible in comparison with the approximation error. Error estimates for RBF approximations of PDE solutions were studied by several authors~\cite{Safdari,RBFLS} in cases when isotropic Gaussian and inverse multiquadric basis functions are used. They demonstrated that the PDE error is governed by the RBF interpolation error
\begin{equation}\label{Interp_error}
\mathcal{E_L} = \mathcal{L}\big(v-\mathcal{J}(v)\big).
\end{equation}
Here $\mathcal{J}(v)$ is the RBF-PUM interpolant defined as
\begin{equation}
\mathcal{J}(v) = \sum_{k=1}^{M}w_k\mathcal{J}(v_k),
\end{equation}
where $\mathcal{J}(v_k)$ is the local RBF interpolant defined as in~\eqref{LocRBFinterp} satisfying the interpolation condition 
\begin{equation}
\mathcal{J}(v_k)(\underline{x}^{k}) = v(\underline{x}^{k}),
\end{equation}
where $\underline{x}^k$ are the $n_k$ local nodes that belong to the patch $\Omega_k$. Therefore, in this section we aim to adapt the existing isotropic interpolation error estimates for the anisotropic Gaussian RBFs.

In particular, the authors of~\cite{Safdari,Rieger,RBFLS} consider two node refinement strategies: (i) refining the node set while having the domain partitioning fixed; and (ii) refining the node set and simultaneously refining the domain partitioning. Strategy (i) leads to an increasing number of local nodes in each patch and, thereby, yields an exponential convergence (for smooth problems). The limitation of this strategy is that after just a few iterations of node refinement the coefficient matrix becomes highly ill-conditioned. Strategy~(ii) leads to an increasing number of patches that contain a similar number of local nodes in every iteration of node refinement. This strategy yields only an algebraic convergence, but the condition number of the coefficient matrix does not increase as rapidly. In this paper \mbox{we are interested in} and adhere to Strategy (ii) because it allows for larger scale simulations. 

Thus, recapitulating the results from~\cite{Safdari,RBFLS} for isotropic RBF-PUM we have
\begin{equation}\label{isotropicBound}
||\mathcal{E_L}||_{L_{\infty}(\Omega)} \leq K \max_{1\leq k\leq M}\sum_{|\beta|\leq |\alpha|}\binom{\alpha}{\beta}||D^{\beta} w_k||_{L_{\infty}(\Omega_k)}  ||D^{\alpha-\beta}\big(v_k-\mathcal{J}(v_k)\big)||_{L_{\infty}(\tilde\Omega_k)},
\end{equation}
where $\tilde \Omega_k = \Omega\, \cap \, \Omega_k$, $v_k$ is the global solution restricted to $\tilde \Omega_k$, $\alpha$ is the  degree of the differential operator $\mathcal{L}$, defined as a multi-index $\alpha = (\alpha_1,\ldots,\alpha_d)$, \mbox{$D^{\alpha}$ is the differentiation} operator of degree~$\alpha$, and $K$ is some constant. The partition of unity weight functions $w_k$ are required to have bounded derivatives up to order $\alpha$, i.e.,
\begin{equation}\label{DiamBound}
||D^{\alpha} w_k||_{L_{\infty}(\Omega_k)} \leq \frac{C_{\alpha}}{H^{|\alpha|}_k},
\end{equation}
where $H_k$ is the diameter of $\Omega_k$. From~\eqref{isotropicBound} using the results of~\cite{Rieger} combined with~\eqref{DiamBound} one can obtain the following error estimate, which indicates an algebraic rate of convergence,
\begin{equation}\label{AlgConv}
||\mathcal{E_L}||_{L_{\infty}(\Omega)} \leq K \max_{1\leq j\leq M} C^{A}_k H_k^{q(n_k)+1-\frac{d}{2}-|\alpha|}||v||_{\mathcal{N}(\tilde\Omega_k)},
\end{equation}
where the norm $||\cdot||_{\mathcal{N}(\tilde\Omega_k)}$ denotes the native space norm (see~\cite{Fasshauer,Rieger}) associated with the chosen RBF, the function
$q(n_k)$ corresponds to the polynomial degree $q$ supported by the local number of points $n_k$, $d$ is the dimensionality of the problem, and $C^{A}_k$ are some constants that depend on the dimensionality $d$, the chosen weight function, the number of local points $n_k$, and the order of the differential operator.

 Introducing anisotropic scaling is equivalent to changing the distance function, in which RBFs are defined. The properties of the functions remain unchanged. Differentiation of an anisotropic RBF  gives rise to the aspect ratio constant coming from the distance defined in $||\cdot||_a$. The same applies to differentiation of the partition of unity weight functions. Hence, we can adapt  \eqref{isotropicBound} in the following way
\begin{align}\label{anisotropicBound}
||\mathcal{E_L}||_{L_{\infty}(\Omega)} \leq K \max_{1\leq j\leq M}\sum_{|\beta|\leq |\alpha|}\binom{\alpha}{\beta} & a^{\beta}||D^{\beta} w_k||_{L_{\infty}(\Omega_k)} \, \times \nonumber \\
&a^{\alpha-\beta}||D^{\alpha-\beta}\big(v_k-\mathcal{J}(v_k)\big)||_{L_{\infty}(\tilde\Omega_k)}.
\end{align}
where
\begin{equation}
a^{\beta} = a_1^{\beta_1}\times\ldots\times a_d^{\beta_d}.
\end{equation}
That is, in total we get the factor $a^{\alpha}$ appearing in the inequality. Thereby, we can again use the results of~\cite{Rieger} to obtain the following error estimate that adapts for the anisotropic RBF-PUM
\begin{equation}\label{anisotropicAlgConv}
||\mathcal{E_L}||_{L_{\infty}(\Omega)} \leq a^{\alpha} K \max_{1\leq j\leq M} C^{A}_k H_k^{q(n_k)+1-\frac{d}{2}-|\alpha|}||v||_{\mathcal{N}(\tilde\Omega_k)}.
\end{equation}
The relation~\eqref{anisotropicAlgConv} goes in line with the results for error bounds for the interpolation by anisotropic RBFs obtained in~\cite{Beatson2010}. Also, it indicates that the approximation by the anisotropic RBF-PUM converges algebraically if the refinement strategy (ii) is employed keeping the norm fixed.

\section{Numerical Experiments}\label{sec: results}

In this section we present results of numerical experiments for the ISMIP-HOM B benchmark test and for  the two-dimensional ice cap. The ISMIP-HOM~B test is a glacier size problem, whereas the two-dimensional ice cap represents a continental size ice sheet.
We illustrate the advantages of the anisotropic radial basis function methods over the standard finite element method (with piecewise linear basis functions and anisotropic unstructured mesh) on the ISMIP-HOM~B test.
In the two-dimensional ice cap simulation, we demonstrate the capability of the anisotropic RBF partition of unity method to solve a continental size ice sheet problem, where the ratio between the ice length and height is around $428:1$. We would like to emphasize here that isotropic RBF methods fail to provide any adequate solution. The reason is the impossibility to select a suitable shape parameter that would define basis functions of the appropriate shape to properly resolve in both directions. Preparing this paper we tested isotropic RBF methods on several different node layouts with various shape parameters. However, in none of the tests we succeeded to obtained a stable solution. Here we do not provide the velocity fields approximated by isotropic RBF methods because they make no sense.
Potentially, one could use the same node set resolution in both $x$- and $z$-directions, but then having, say, $35$ nodes in the vertical direction with the maximal ice thickness of 3.5 km would require having 15000 nodes in the horizontal direction to provide the same node set resolution, if the ice length is 1500 km. Such a resolution would lead to a high computational cost of the algorithm. Therefore, we are certain that anisotropic RBF methods are much more suitable for approximation of the dynamics of continental size ice sheets.


In order to construct a set of computational nodes we use (i) a background Cartesian grid and (ii) a set of Halton nodes. The Cartesian grid is characterized by $N_x$ and $N_z$ nodes in the $x$- and $z$-directions with resolutions $h_x$ and $h_z$, respectively, which are distributed over the rectangle $[0,L_x]\times[0,L_z]$, where $L_x$ is the length of the domain and $L_z$ is the maximum thickness. The quasi-random Halton nodes are also selected to fill up the rectangle $[0,L_x]\times[0,L_z]$ with $N_x\times N_z$ nodes. To define the resolution of the Halton node set we simply introduce the notation (likewise for the Cartesian grid) $h_x=L_x/(N_x-1)$ and $h_z=L_z/(N_z-1)$. After defining nodes for the rectangle, the nodes which fall outside the domain profile are discarded, and the ones which fall inside are kept and augmented by the boundary nodes, thereby resulting in a total of $N$ computational nodes.

The Cartesian background grid is very convenient and easy to construct. However, it is not always the most suitable node layout for the RBF approximation~\cite{Fornberg2010} but allows to select the shape parameter $\varepsilon$ relatively simply since the distance in both $x$- and $z$-directions is uniform in the $||\cdot||_a$-norm and does not vary. In contrast, the Halton nodes might be a bit more challenging to select a proper shaper parameter for. Although, the Halton nodes have a good space filling property and are, thus, maybe more suitable for large simulations. 


So far no analytical procedures of finding the shape parameter have been developed. Therefore, we seek the shape parameter empirically in the following form
\begin{equation}\label{ShPar}
\varepsilon = \frac{C}{h},
\end{equation}
where $C$ is a constant and $h$ is the internodal distance defined in the $a$-norm as
\begin{equation}
  h=\sqrt{h_x^2+a^2h_z^2}.
\end{equation}
and the aspect ratio $a$ is defined as the ratio between the spacings in the $x$- and $z$-directions.
\begin{equation}
a = \frac{h_x}{h_z}.\label{eq: aspectR}
\end{equation} 
We find that $C=0.5$ is a fairly good choice for the ISMIP-HOM B benchmark test and that $C=0.125$ is a suitable choice for the ice cap test case. In order to get an intuition on how to choose the value of the constant $C$ we adopt a strategy of controlling the condition number of the interpolation matrix $A$ (in case of partition of unity this is the global interpolation matrix reassembled from the local interpolation matrices with the partition of unity weights). We maintain the condition number around $10^{16}$, i.e., around near ill-conditioning, as it is a well-known fact that the RBF approximation achieves its best accuracy for small values of $\varepsilon>0$~\cite{Larsson7,Fornberg3,Larsson2003}, which leads to a large condition number of the matrix. We repeat this procedure several times on coarse node sets and record the condition numbers of the resulting matrices. Then we utilize a least squares fit to determine the value of~$C$ and reuse this value for finer node sets. To estimate the condition number we use the \mbox{1-norm} condition estimator of Hager~\cite{Hager1984} and a block-oriented generalization of Hager's estimator~\cite{Higham2000} provided by the function \texttt{condest} in MATLAB. Some other strategies and insights on how to determine a proper value of the shape parameter can be found in~\cite{Ahlkrona,Milovanovic}.

\subsection{ISMIP-HOM B}
\label{sec:ISMIP-B}

We compute the velocity field of an ice slab that is grounded on a bedrock that has sinusoidal shape. The length of the slab is $L=10$~km.
The surface is defined as a function of the $x$-coordinate
\begin{equation}
	h_s(x)=-x\tan{\alpha}, \quad x\in[0,L],\label{eq:surf ele}
\end{equation}
and has a slope $\alpha=0.5^\circ$.
The bedrock is defined as
\begin{equation}
	h_b(x) = h_s(x)-1000+500\sin\left(\frac{2\pi x}{L}\right), \quad x\in[0,L].
\end{equation}

\begin{figure}[h]
\label{ISMIPrbffem}
  \center 
  \includegraphics[width=1\textwidth, trim={80 0 0 0},clip]{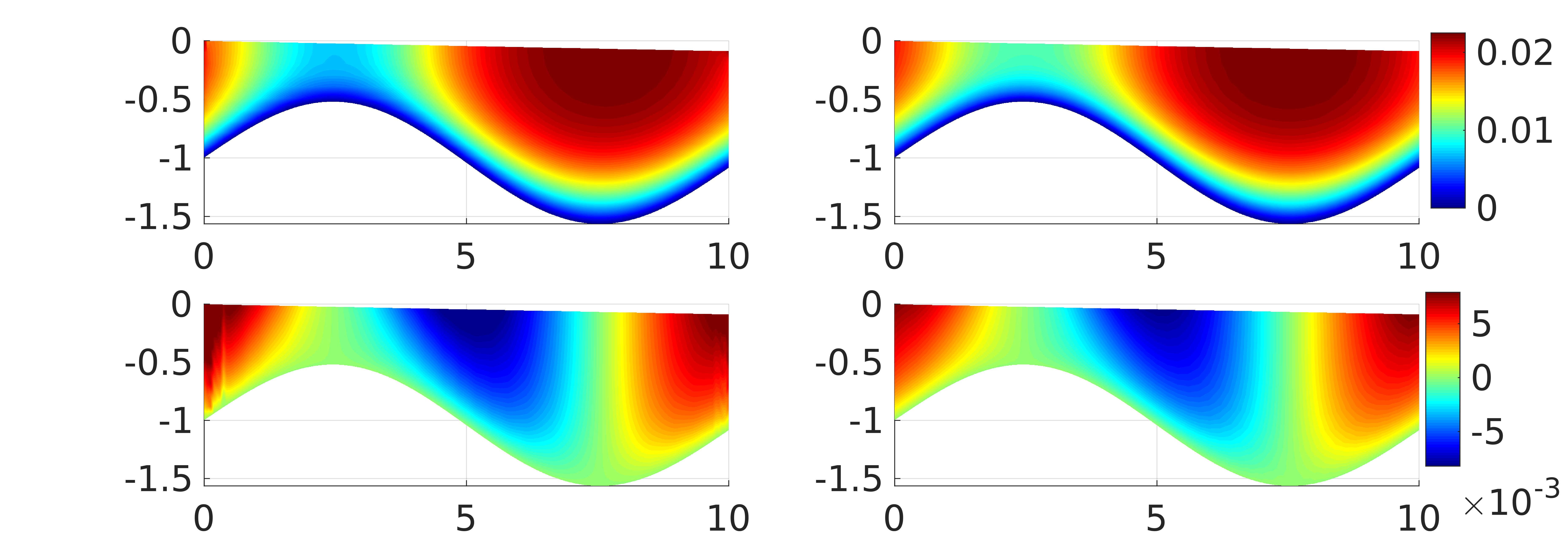} 
  \caption{\emph{Top:} The horizontal velocity for the ISMIP-HOM B test computed by the global RBF method (\emph{left}) and FEM (\emph{right}). \emph{Bottom:} The vertical velocity computed by the global RBF method (\emph{left}) and FEM (\emph{right}). All velocities are measured in  km/a.} 
\end{figure}

We assign periodical boundary conditions at the lateral boundaries. Also, we assume that the ice is grounded on the bedrock, meaning no slip boundary condition on the bottom, and the ice surface is free of stress, which leads to a free surface problem.

The horizontal and vertical velocities, obtained by the global RBF method and FEM on node sets with similar resolution, are presented in Figure~\ref{ISMIPrbffem}. The two left panels represent the RBF solutions, while the right panels represent the FEM solutions, which we use for comparison. The horizontal and the vertical velocities are well resolved by both methods.

%

%
\begin{table}[h!] 
\begin{center}
\caption{Maximum absolute errors (m/a) in the horizontal velocity and CPU times for several different numbers of degrees of freedom for FEM, global RBF method, and RBF-PUM for the ISMIP-HOM B benchmark test.}
\label{TabRBFFEM}
\resizebox{\textwidth}{!}{\begin{tabular}{c r r | c r r | c r r}
\hline\hline 
  \multicolumn{3}{ c |}{FEM} & \multicolumn{3}{ c |}{RBF} & \multicolumn{3}{ c }{RBF-PUM} \\   \hline\hline 
 Error  & \multicolumn{1}{ c }{$N$} &  Time (s) &  Error & \multicolumn{1}{ c }{$N$} &  Time (s) &  Error & \multicolumn{1}{ c }{$N$} & Time (s)\\  \hline
 $0.3353$  & $243$  & $2.2782$ & $0.2535$  & $306$ &  $0.2388$ & $0.3467$  & $291$ &  $0.2479$ \\
 $0.1528$  & $805$  & $8.0495$ & $0.1574$  & $802$ &  $1.3701$ & $0.1365$  & $771$ &  $0.7685$ \\
 $0.0492$  & $2889$  & $34.5913$ & $0.0432$  & $1992$ &  $13.7930$ & $0.0372$  & $1925$ &  $4.1405$  \\
\hline\hline
\end{tabular}}
\end{center}
\end{table}

We compare three methods: the global RBF method, RBF-PUM, and FEM in terms of accuracy and computational efficiency. We use a finite element solution on a fine grid with $10897$ degrees of freedom as our reference solution. The maximum absolute errors, run times, and numbers of degrees of freedom are presented in Table~\ref{TabRBFFEM}. We observe that the RBF methods give better accuracy in much shorter time. RBF-PUM is more than $8$ times faster than the standard FEM. Additionally, thanks to the high accuracy, the RBF methods need fewer nodes to reach a similar error tolerance level. This is crucial for large simulations since less storage will be required. The reason that the FEM implementation runs slower is the necessity of matrix reassembly within the nonlinear iteration. It has been shown~\cite{Ahlkrona2} that this procedure may severely dominate the total computational time.


Thereby, the numerical results for the ISMIP-HOM B demonstrate that RBF-PUM is the most efficient out of the three methods for this type of problem. Inspired by this, we continue investigating the properties of the anisotropic RBF-PUM. However, in the remaining part of the paper we will not provide comparisons with the global RBF method and FEM, because the global RBF method becomes too computationally intensive, while a rigorous implementation of FEM for the two-dimensional ice cap goes beyond the scope of this paper.

\subsection{The Two-Dimensional Ice Cap}
\label{sec: 2D ice}

This experiment is inspired by the EISMINT benchmark test \cite{huybrechts_eismint_1996}, where numerical methods are intercompared on a continental size ice sheet. There are several challenges in simulating ice dynamics with RBF methods even at steady states, for instance, a large aspect ratio of the computational domain or a dramatic variation of the ice thickness from the margins to the centre of the cap. These factors would require a delicate adjustment of the shape parameter in the isotropic RBF case, but it is much simpler for the anisotropic RBFs since the issues with the aspect ratio can be resolved by introducing the distance in the $||\cdot||_a$-norm. Thus, in this test problem we would like to  illustrate the clear advantages of anisotropic RBF methods over the standard isotropic RBF methods.


The computational domain follows the so-called Bueler profile \cite{greve2009dynamics,bueler2005exact}, which is a two-dimensional flow-line model spanning from $x=0$~km to $x=1500$~km. At steady state, the ice cap covers the area from $x=300$~km to $x=1200$~km and the ``ice free'' area is at $x\in[0,300]$ and $x\in[1200,1500]$, (see Figure~\ref{fig:EISMINT_PUM}). Generally, in order to make it possible for the \mbox{margin of the} ice cap to extend, a thin layer ($10$~m) of ice is assumed on the ``ice free'' area.

The physical parameters of this experiment are given in Table \ref{tab: physical param}. There is no sliding on the bedrock and no velocity on the lateral boundaries. The top surface is free of stress by neglecting the atmospheric pressure. In order to accurately capture the surface gradient Chebyshev points are used to generate the upper boundary discretisation on the ice cap part. The ``ice free'' regions and the bottom boundaries are discretised using equidistant nodes.



\subsubsection{The Cartesian Nodes}

The anisotropic RBF provides the flexibility to choose any combination of $N_x$ and $N_z$. However, we need to keep in mind that the expression \eqref{eq: aspectR} for the aspect ratio relates $N_x$ and $N_z$ as 
\begin{equation}
  \frac{N_x}{N_z}=\frac{L_x}{aL_z}.
\end{equation}
Here the aspect ratio can be considered as a scaling of the $z$-direction. Ideally, we prefer to have the final geometry after the scaling close to a square. That is, we need to select $N_x$ and $N_z$ with some special care. Therefore, in this application we require the numbers $N_x$ and $N_z$ to be of the same order of magnitude.

Also, note that when combining the background and boundary nodes we need to account for nearly indistinct nodes, i.e., nodes positioned in a tiny neighbourhood of each other, since this can lead to a highly ill-conditioned or singular coefficient matrix. In order to avoid the ill-conditioning we remove the internal nodes whose distance to the boundary is shorter than $h/4$. 

 \begin{figure}[!h]
  \center 
  \includegraphics[width=0.8\textwidth]{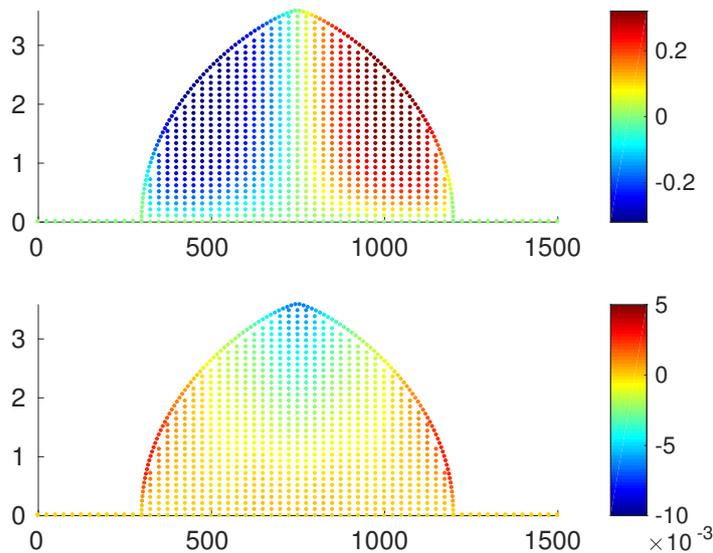} 
  \caption{\emph{Top:} The horizontal velocity (in km/a) of the ice cap  computed on Cartesian nodes by anisotropic RBF-PUM. \emph{Bottom:} The respective vertical velocity (in km/a). } \label{EISMINT_Cartesian}
\end{figure}

The result of a numerical simulation with anisotropic RBF-PUM is presented in Figure~\ref{EISMINT_Cartesian}. The horizontal velocity is displayed in the upper panel and the vertical velocity is displayed in the lower panel. The unit of the velocities is kilometre per annum (km/a). The background grid is characterised by $N_x=60$ and $N_z=35$, yielding $N=1014$ computational nodes in the domain. The reader is referred to \cite{huybrechts_eismint_1996} to ascertain that the presented solutions are qualitatively correct.

\subsubsection{The Halton Nodes}
\label{sec: halton}

Another valuable property of RBF methods is their mesh-free nature, which makes them suitable for problems set in domains with complex geometries. Moreover, having the mesh-free property we can cluster computational nodes and achieve a higher resolution in the regions where it is required, e.g., along the margins. In order to demonstrate the capability of anisotropic RBF methods to work on unstructured node sets we select the quasi-random Halton node layout. The scaling of the anisotropic basis functions is based on aspect ratio as a function of the mesh resolutions as in relation \eqref{eq: aspectR}.



\begin{figure}[!h]
  \center 
  \includegraphics[width=0.8\textwidth]{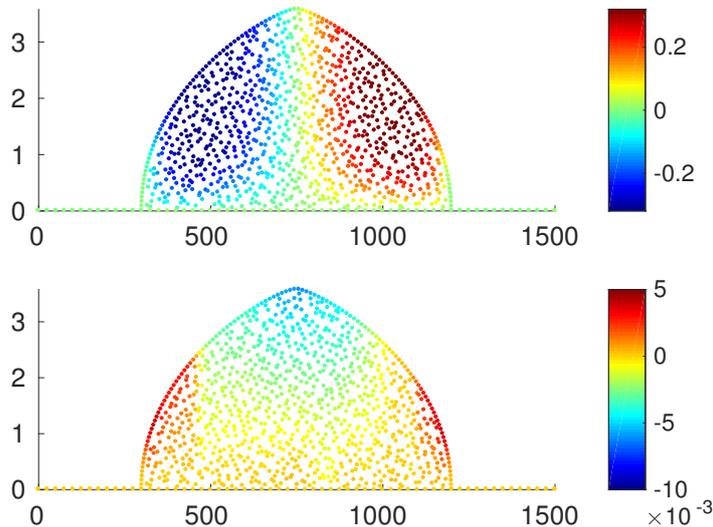} 
  \caption{\emph{Top:} The horizontal velocity (in km/a) of the ice cap  computed on the Halton nodes by anisotropic RBF-PUM. \emph{Bottom:} The respective vertical velocity (in km/a). } \label{EISMINT_Halton}
\end{figure}

We generate and scatter the Halton nodes over the computational domain according to the strategy given in the beginning of Section~\ref{sec: results}. After we augment them with the boundary nodes and remove the indistinct nodes to avoid singularities in the coefficient matrix. In total we obtain 1063 nodes in the computational domain, which are used to construct the approxiamtion. In Figure~\ref{EISMINT_Halton} we display the solutions by anisotropic RBF-PUM. Both horizontal and vertical velocities are properly resolved and the results go in line with the results obtained on the Cartesian nodes.

\subsubsection{Convergence Tests}




In this section, we study convergence of the anisotropic RBF-PUM for both Cartesian and quasi-random nodes. The reference solution is computed by the anisotropic RBF-PUM on $23943$ Cartesian nodes.
To generate a convergence experiment we repeatedly solve for the horizontal velocity gradually refining the node set and measure the root mean square error against the reference solution. The inequality~\eqref{anisotropicAlgConv} holds for the root mean square error since 
\begin{equation}
\frac{1}{\sqrt{N}}||\mathcal{E_L}||_{L_{2}(\Omega)} \leq ||\mathcal{E_L}||_{L_{\infty}(\Omega)}.
\end{equation}

For the convergence experiments we used a mode with approximately 150 nodes per internal patch. However, the patches, which are close to the ice margins, contain just around 20-25 nodes. Therefore, the overall convergence rate was dictated by the order of the local error convergence in those patches. Having the number of local nodes $n_k=21$ means that the highest polynomial degree that can be supported in two-dimensional space is $q(n_k)=5$. That is, taking into account that the problem dimensionality $d=2$ and the order of the differential operator $|\alpha|=2$, we conclude that the anisotropic RBF-PUM approximation should converge with order 3.

\begin{figure}[!ht]
  \center 
    \includegraphics[width=0.85\textwidth]{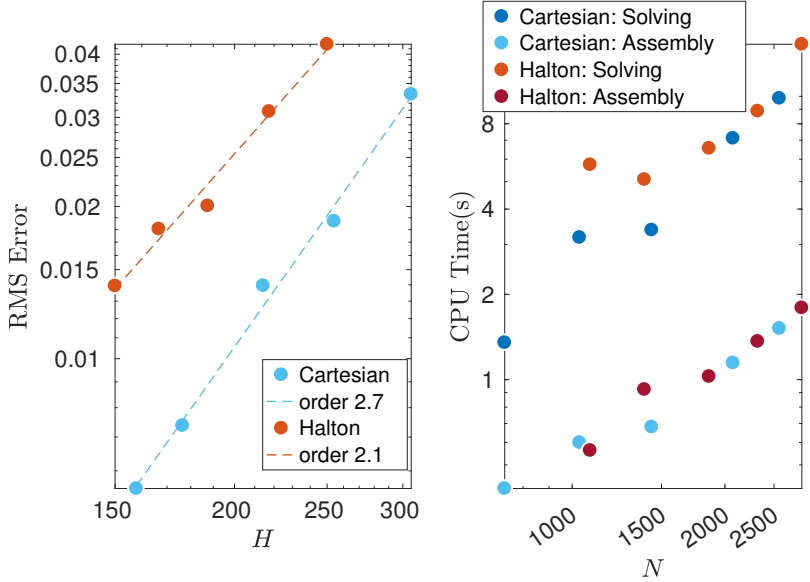} 
  \caption{\emph{Left:} The convergence of the approximation of the horizontal velocity obtained by the anisotropic RBF-PUM with respect to the patch sizes $H$ on the Cartesian nodes (blue) and Halton nodes (red) in the log-log scale. \emph{Right:}  The computational time of the matrix assembly and solving procedure with respect to the total number of computational nodes in the log-log scale.} 
  \label{fig: convergence}
\end{figure}

As we can see in the left panel of Figure~\ref{fig: convergence}, the rate of convergence on the Cartesian nodes with respect to the patch radius is around~$2.7$, which goes in line with estimate~\eqref{anisotropicAlgConv} as well as with the results in~\cite{Ahlkrona}. The rate of convergence on the Halton nodes is slightly lower and is around~$2.1$, which makes it require about $\sqrt{2}$ times more  nodes to reach a similar quality of the solution. Roughly speaking, with a patch resolution of 150~km we get the horizontal velocity resolved up to 6~m on the Cartesian nodes and up to $14$~m on the Halton nodes.

\begin{figure}[!ht]
  \center 
    \includegraphics[width=1\textwidth]{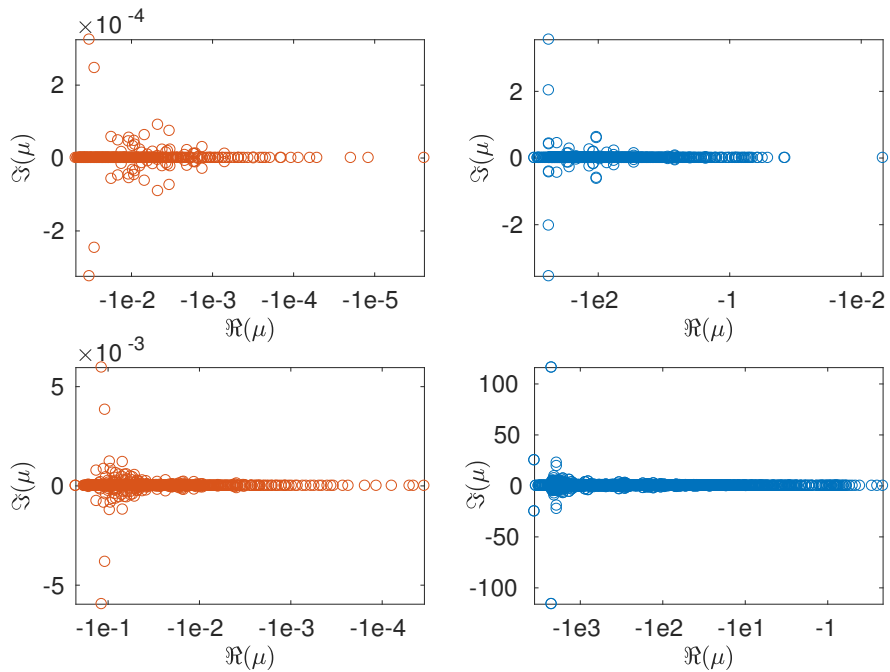} 
  \caption{\emph{Top:} Spectra of the differential operator $\mathcal{L}$ discretized by isotropic RBFs (\emph{left}) and anisotropic (\emph{right}) RBFs on a coarse node set with $736$ Cartesian nodes. \emph{Bottom:} Their counterparts for $\mathcal{L}$ discretized on a fine node set with $2562$ Cartesian nodes.} 
  \label{fig: eigenvalues}
\end{figure}

In general, we could employ adapted node layouts and patches with the size matching the domain profile to contain the same number of nodes regardless of their position to relax the constraint on the convergence rate, but we leave it to be implemented elsewhere. Also, in general, our experience shows that one should not expect very high convergence rates for this type of problems. The considered problem is nonlinear with the viscosity inversely proportional to the velocity, meaning that it approaches infinity for low velocities, which introduces a singularity in the formulation. In order to avoid the singularities and proceed with the numerical solution we prevent the viscosity of growing above $10^{10}$ by setting a threshold. However, this turns the viscosity into a non-smooth function, which usually does not allow for getting optimal convergence rates.

We present the computational cost of assembling the RBF matrices and solving the nonlinear system in terms of  CPU time in Figure~\ref{fig: convergence}. The CPU time is averaged over five independent runs for each resolution to eliminate system noise. For both Cartesian and Halton nodes the time for the matrix assembly is a linear function of the total number of nodes since the partition of unity method provides a sparse structure of the system and the number of nodes per partition remains fixed such that the bandwidth of the sparse matrix does not grow with the total number of nodes. In the solving procedure, we use the default backslash operator in MATLAB. The cost of solving the nonlinear system increases linearly or perhaps slightly super-linearly with the total number of computational nodes thanks to the sparse structure of the coefficient matrix (opposed to $\mathcal{O}(N^3)$ for dense systems).  

Additionally, we look at the spectra of the linearized Blatter--Pattyn operator $\mathcal{L}$ discretized on two Cartesian node sets (coarser and finer) by both isotropic and anisotropic RBF-PUM and present them in Figure~\ref{fig: eigenvalues}. We use the same RBF for both methods. The only difference is that the RBF is evaluated with respect to the Euclidian norm for the isotropic case and with respect to the $||\cdot||_a$-norm for the anisotropic case. For the operator $\mathcal{L}$ we would like to see all eigenvalues in the negative half plane, close to the real axis. We observe that both isotropic and anisotropic approximations give spectra, which are rather aligned along the real axis. However, we also notice that the isotropic RBF-PUM yields a spectrum, which is more clustered towards the origin. In fact, using MATLAB's \texttt{eig} function for computing matrix eigenvalues, we find that there is at least one eigenvalue with $\Re(\mu)=0$. That is, the discrete operator is numerically singular. Therefore, the isotropic RBF-PUM is not capable to provide any sensible approximation. In contrast, the anisotropic RBF-PUM yields eigenvalues, which are distinct from the origin, and the ratio $\max\big(|\Re(\mu)|\big)/\min\big(|\Re(\mu)|\big)$decreases under the node set refinement and is $2\cdot10^5$ for a set of 736 Cartesian nodes and $1.7\cdot10^4$ for a set of 2562 Cartesian nodes.

\section{Conclusions}\label{sec:conclusion}


The goal of this work was to extend the approach in~\cite{Ahlkrona} for domains with extreme length/thickness ratios to enable approximating the velocities of continental size ice sheets. The main obstacle was the incapability of standard (isotropic) RBF methods to cope with different scales of ice sheet modelling. In order to overcome this issue we developed and implemented an anisotropic RBF partition of unity method that accounts for the aspect ratio of  the domain discretisation and incorporates the ratio into the basis function structure through the updated distance function. This allows us to modify the basis functions in such a way to best fit the domain geometry.

The anisotropic RBF-PUM was tested on a glacier size benchmark test ISMIP-HOM~B and on a synthetic continental size ice sheet whose geometry was inspired by the EISMINT benchmark. The method was compared with the global RBF method and FEM for the glacier size test case and demonstrated a much greater efficiency than its counterparts. 

To recover the velocity field of the ice cap we used Cartesian and Halton nodes. The Cartesian nodes were found to give a better accuracy than the Halton nodes. However, the difference in the results was not critical. In fact, our main idea of implementing the Halton nodes was to demonstrate that the anisotropic RBF-PUM is suitable for both structured and unstructured node sets. 

Additionally, we adapted the error estimates from~\cite{Safdari,RBFLS} for the anisot-ropic RBF-PUM. We tested the convergence of the solutions obtained by the anisotropic RBF-PUM on the two node layouts. The convergence rates on the Halton and Cartesian nodes were~2.1 and 2.7, respectively. These experimental results go in line with the analytical estimates. In order to select the shape parameter value for the refined node sets we developed a strategy bases on controlling the condition number of the interpolation matrix. 


We showed that the anisotropic RBF-PUM is suitable for different scales of ice sheet modeling. The method is more efficient and accurate than the finite element method in the investigated cases. It becomes even more important for applications with very high aspect ratios since special treatments are required for FEM which may lower the accuracy or increase the computational complexity. By using the anisotropic RBF-PUM, all the good features of the isotropic RBF-PUM are preserved without extra cost, while new crucial features are gained.


This paper together with~\cite{Ahlkrona} provides a base for an RBF framework for ice sheet simulations that comprises frozen basal and discontinuous partial slip basal conditions, and moving surface. The method is suitable for various geometries, even with extreme length/thickness ratios, and capable for approximating the ice velocity on both structured and unstructured sets of computational nodes, which is an important property for developing adaptive strategies and taking it further to, say, a grounding line problem~\cite{pattyn_results_2012}.

\newpage
\section*{\large{Acknowledgements}}

\footnotesize{The authors would like to thank Josefin Ahlkrona, Elisabeth Larsson, Jeremy Levesley, Igor Tominec, Lina von Sydow, and Grady Wright for fruitful discussions. Moreover, the authors are grateful to Elisabeth Larsson and Lina von Sydow for proofreading the manuscript. 
Cheng Gong was supported by the Swedish Research Council FORMAS grant 2013-1600.}






\bibliographystyle{unsrt}
{\footnotesize
\bibliography{reference}}

\end{document}

%% file: Util.tex
\usepackage{listings}

\usepackage[latin1]{inputenc}
\usepackage[english]{babel}
\usepackage{ifpdf}
\usepackage{amsmath}
\usepackage{amssymb}
\usepackage{amsthm}
\ifpdf
  \usepackage[pdftex]{graphicx}
  \usepackage{epstopdf}
\else
  \usepackage[dvips]{graphicx}
\fi

\usepackage{float}
\usepackage{tikz}
\usepackage{color}
\usepackage[pdftitle={Manuscript},
  pdfauthor={Cheng Gong},
  pdffitwindow=true,
  breaklinks=true,
  colorlinks=true,
  urlcolor=blue,
  linkcolor=black,
  citecolor=black,
  anchorcolor=black]{hyperref}

\definecolor{hellgelb}{rgb}{1,1,0.85}
\definecolor{colKeys}{rgb}{0,0,1}
\definecolor{colIdentifier}{rgb}{0,0,0}
\definecolor{colComments}{rgb}{0.2,0.7,0}
\definecolor{colString}{rgb}{0.788,0,0.961}
\newcommand{\pdd}[2]{\frac{\partial #1}{\partial #2}}

